# Decision Making and Productivity Measurement

Dariush Khezrimotlagh
Department of Statistics, Miami University, Ohio, USA

**Abstract**

In this article, the concepts of technical efficiency, efficiency, effectiveness and productivity are illustrated. It is discussed that when firms are not homogenous, the situation is the same as when each factor has a different unit of measurement from one firm to another, and therefore, no meaningful discrimination can be expressed, unless a set of known weights are introduced to standardize data. A linear programming DEA model is used when a set of known weights are given to calculate the technical efficiency and efficiency of a set of homogenous DMUs with multiple input factors and output factors. A numerical example is also provided.

## 1.0. Introduction

A non-parametric technique to estimate the technical efficiency and the efficiency of a set of homogenous firms was proposed by Farrell (1957). His method estimates the production function non-parametrically, which was similarly suggested by Debreu (1951) and Koopmans (1951). Førsund and Hjalmarsson (1974) illustrated the notion of efficiency in the macro, the industry and the micro levels, and clearly displayed and demonstrated the differences between the production frontier and efficiency. Soon after, Charnes *et al*. (1978) proposed a mathematical construction and a linear programming model to introduce technically efficient firms with multiple input factors and multiple output factors. They called the mathematical construction '*Data Envelopment Analysis'* (DEA) and the model '*Charnes, Cooper and Rhodes*' (CCR). CCR generates a PPS based upon a set of available homogenous firms, and non-parametrically and linearly estimates the production function; thus, the firms which lie on the frontier of that PPS are called technically efficient. CCR fails to measure the technical inefficiency completely, and only calculate the output view of technical efficiency (that is, increasing the values of output factors without measuring the excess of input factors) or calculate only the input view of technical efficiency (that is, decreasing the values of input factors without measuring the shortage of output factors). For this reason, Färe and Lovell (1978) noted on Farrell's measurement of technical efficiency and CCR, and proposed a *Russell measure*, to simultaneously deal with both input and output views of technical efficiency. Their proposed model is a non-linear programming, and difficult to solve; thus, Pastor *et al*. (1999) proposed *Enhanced Russell Measure* (ERM) to measure technical efficiency of firms and avoid computational and interpretative difficulties with the Russell measure.

On the other hand, Charnes *et al*. (1985) proposed an *Additive model* (ADD) to remove the shortcomings of CCR to measure the technical inefficiency of firms. Nonetheless, ADD is also not a perfect model to measure the technical inefficiency. Therefore, Tone (1997, 2001) proposed a *Slack-Based Measure* (SBM) model to measure technical inefficiency of firms. He proved that (1) ERM and SBM are equivalent in that the lambda's values that are optimal for one are also optimal for the other, (2) the SBM measurement corresponds to the mean proportional rate of input factors' reduction and the mean proportional rate of output factors' expansion, (3) the SBM measurement is monotone, decreasing in each input and output slack, and (4) it is invariant with respect to the unit of measurement of each input and output item.

All proposed models, such as CCR, ADD and SBM are provided to measure the technical inefficiency which neither provides a ranking and benchmarking tool, nor introduces the relative efficiency scores for firms. Note that there are a lot of proposed models based on CCR in the literature of operations research since 1978 which have the same (or more) mentioned shortcomings to focus on technical efficiency (doing the job right) instead of efficiency (doing the job well), to discriminate a set of homogenous firms. While these studies not be further upon here, but readers can examine several of these questionable studies, and criticize their lack of discrimination and ranking tools as simple exercises.

Sexton *et al*. (1986) wisely noted the shortcomings of CCR and stated that DEA cannot be used to analyze or comment on a firm's (price) efficiency and a firm can be technically efficient, but (price) inefficient. Thus, they proposed a *cross efficiency* model, which was supposed to measure the score that a particular firm receives when it is rated by another firm. Nonetheless, the cross-efficiency score for a firm is also not a relative score for that firm; it is an average value of the relative scores of that firm, according to some specified sets of weights, which neither should be used to rank firms, nor is relatively meaningful. The average values of the relative scores of firms like the maximum (minimum, first quartile, and so on) values of the relative scores, are not relatively meaningful and should not be suggested as the relative scores of firms.

In order to decrease the above shortcomings, Khezrimotlagh *et al*. (2013) proposed a method (called $\varepsilon$-KAM) to bridge the gaps between technical efficiency and efficiency. The score of KAM is different from the scores of other models in the literature, and can be used as a fair judgment tool for ranking and benchmarking firms. When the value of epsilon is 0, the results of 0-KAM identify the firms which are technically efficient and technically inefficient, and should not be used to rank or benchmark firms. As the value of epsilon increases, the results of $\varepsilon$-KAM can be used to rank and benchmark firms, according to the value of epsilon and the introduced assumptions for weights/prices of input and output factors.

The technical efficiency (as well as production function) depends upon the way of introducing the practical points and the efficiency depends upon the weights/worth/prices of the factors, and at least require the combination of the radiate, the convexity and the wholly dominant approaches to linearly estimate the efficiency scores of a set of homogenous firms. In the next sections, after introducing several phrases, $\varepsilon$-KAM is used to measure the efficiency of firms with the least requirements to the radiate, the convexity and the wholly dominant approaches.

## 2.0. Introduced phrases

The purpose of discrimination for a set of homogenous firms (factories, organization, divisions and so on), in which each firm has multiple input factors and multiple output factors, is to find the firms which have lesser values of the input factors and greater values of the output factors. For such an aim, the ratio of a linear combination of the output factors to a linear combination of the input factors is introduced as the measure for the purpose of discrimination. Therefore, the weight/price/worth of each factor is required to introduce the linear combination of the factors and allow measuring the performances of the firms. Since the firms are homogenous, the weight/price/worth of each factor should not be varied from one firm to another, unless all differences known, and multiplied to the values of the factors before the evaluation.

The concepts of ‘doing the job right’, ‘doing the job well’, ‘doing the useful job’ and so on, were similarly introduced in the literature of economics and operations research, but differently interpreted with several words and phrases, such as, *efficiency*, *technical efficiency*, *price efficiency*, *productive efficiency*, *relative efficiency*, *economic efficiency*, *allocative efficiency*, *overall efficiency*, *productivity*, and so on.

Philosophically, we should avoid using several phrases for a concept, and should cautiously clarify whether the relationship between that concept and what we express is meaningful. For instance, the word ‘efficiency’ is commonly used instead of the phrases ‘technical efficiency’ and ‘relative efficiency’. If these phrases illustrate the same concept, the word ‘efficiency’ should be enough to mention that concept and any further terminology is redundant and misleading.

The same criticism can be illustrated for the phrases ‘price efficiency’, ‘overall efficiency’, ‘productive efficiency’ and so on.

According to the Cambridge English dictionary, the word ‘efficiency’ means “the condition or fact of producing the results you want without waste, or a particular way in which this is done”, the phrase ‘technical efficiency’ means “a situation in which a company or a particular machine produces the largest possible number of goods with the time, materials, labor, etc. that are available”, and the word ‘relative’ means “as judged or measured in comparison with something else”. Even from the literal definition, we can see the terms are clearly different in meaning. Therefore, it is vital to review the meaning of these phrases in the literature of economics, engineering and operations research, and reintroduce them with industry-wide accuracy and understanding.

## 2.1. The technical efficiency measurement

Suppose that there are several homogenous firms which each firm uses a set of input factors to produce a set of output factors. In the literature of *the production theory*, a *Production Possibility Set* (PPS) is a set of all possible situations which a set of output factors can be produced from a set of input factors. A *production function* (*production frontier*) is also a function that gives the maximum possible values of the output factors from the values of the input factors. The points on the production function are called the *technically efficient* points, and this definition is matched to the literal definition as well. None of the coordinates of the technically efficient points can be improved without worsening another coordinate, that is, none of the values of the input (output) factors can be decreased (increased) without increasing the value of another input factor or decreasing the value of another output factor. A point which is in the PPS, and does not lie on the production function is called *technically inefficient*. When a point is technically inefficient, at least one of its input or output factors can be improved to reach the production function in order to be technically efficient. In other words, a technical inefficient point is dominated by one technical efficient point at least.

The following figure depicts the production function for a set of six firms, labeled A-F, which each firm has two input factors to produce a single constant output factor, as well as the related PPS and the technically ef-ficient and inefficient firms.

The horizontal axis in the figure represents the values of the first input factor per unit of the output factor and the vertical axis represents the values of the second input factor per unit of the output factor. The curve $SS'$ in the figure is called the production function and the above area of the curve is the related PPS.

The firms A-D which lie on the production function are technically efficient and the firms E and F, which are inside the PPS, are technically inefficient. Firms A-C wholly dominate E; for instance, E and C used the same value of the first input factor, but E used the greater value of the second input factor, and for this reason E is technically inefficient. In other words, the points which are inside the PPS can be compared with at least one point which is on the frontier of the PPS in the figure.

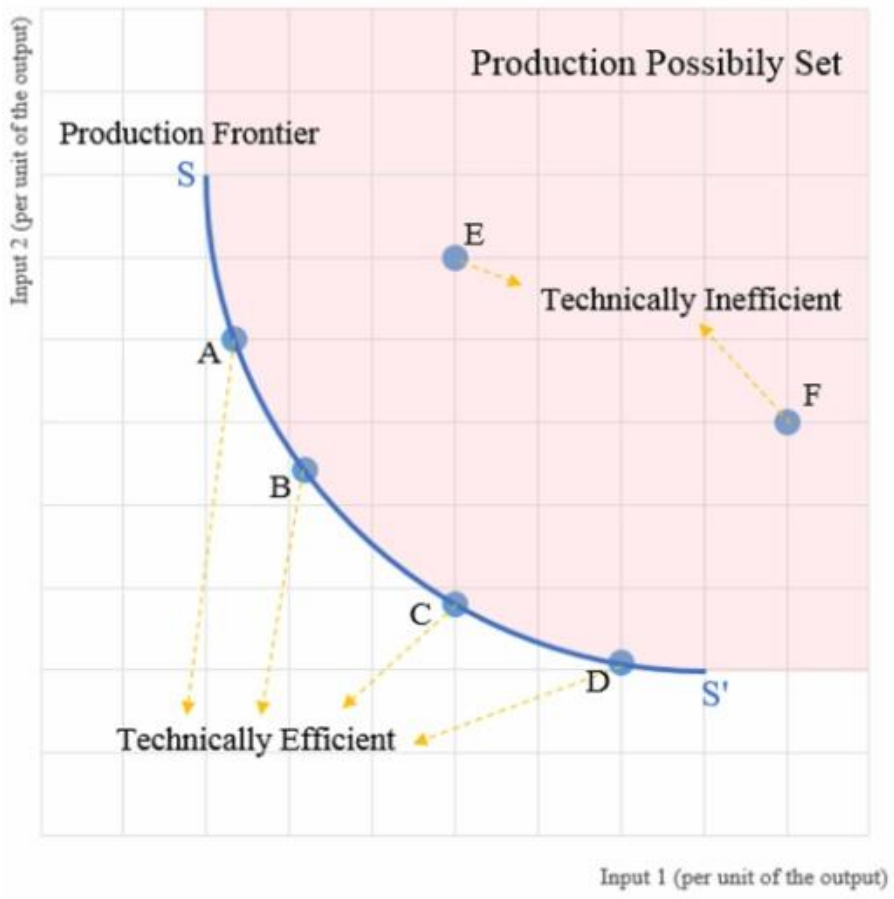

Figure 1: A PPS of a set of six firms.

The PPS is the practical region; the production function is the frontier of the feasible area, and the technically efficient points are the points which have done the job right. In other words, the concept of doing the job right is the same as the concept of technical efficiency, and the word ‘technical’ refers to the used technology (approach) to introduce the practical points.

When a firm has done the job right, it means that the firm has produced the maximum possible values of the output factors from a set of the input factors. In addition, the concept of doing the job right depends upon the introduced approach to generate the practical points, and this is the same as the concept of technical efficiency which depends upon the use of technology to define the production function.

The concept of doing the job right is only a necessary condition to discriminate between firms, and is not enough to introduce the firms which have done the job well. For instance, both firms A and D have done the job right in the above figure; however, if it is supposed that they have done the job well at the same time, a paradox is generated. If a firm lies on the production function, it does not logically say that the firm has done the job well. It is possible that a firm which does not lie on the production function has done the job better than the firm which lies on the production function. In other words, the discrimination between firms based on the production function only, (even if the production function is exactly available), is not valid. The important pros of technical efficiency are to estimate the production function and find the firms which can be candidates for the concept of doing the job well, (without introducing the firm which has done the job well). A firm which has done the job well in comparison with all other firms is technically efficient and lies on the production function, but the points on the production function, which are technically efficient, have not necessarily done the job well.

### 2.2. Efficiency measurement

Let's suppose that the line $TT'$ has the same slope as the ratio of the prices/weights/worth of the two input factors, as depicted in Figure 2. From the figure, D has done the job well in comparison with all other firms. Thus, the firms can be arranged from the highest rank to the lowest rank, given by D, C, B, A, F and E, respectively. Point P is not practical (according to the PPS), but has the same worth as D's performance, and allows discrimination between E and D. E can increase the value of the first input factor and decrease the value of the second output factor to perform as well as D within the PPS.

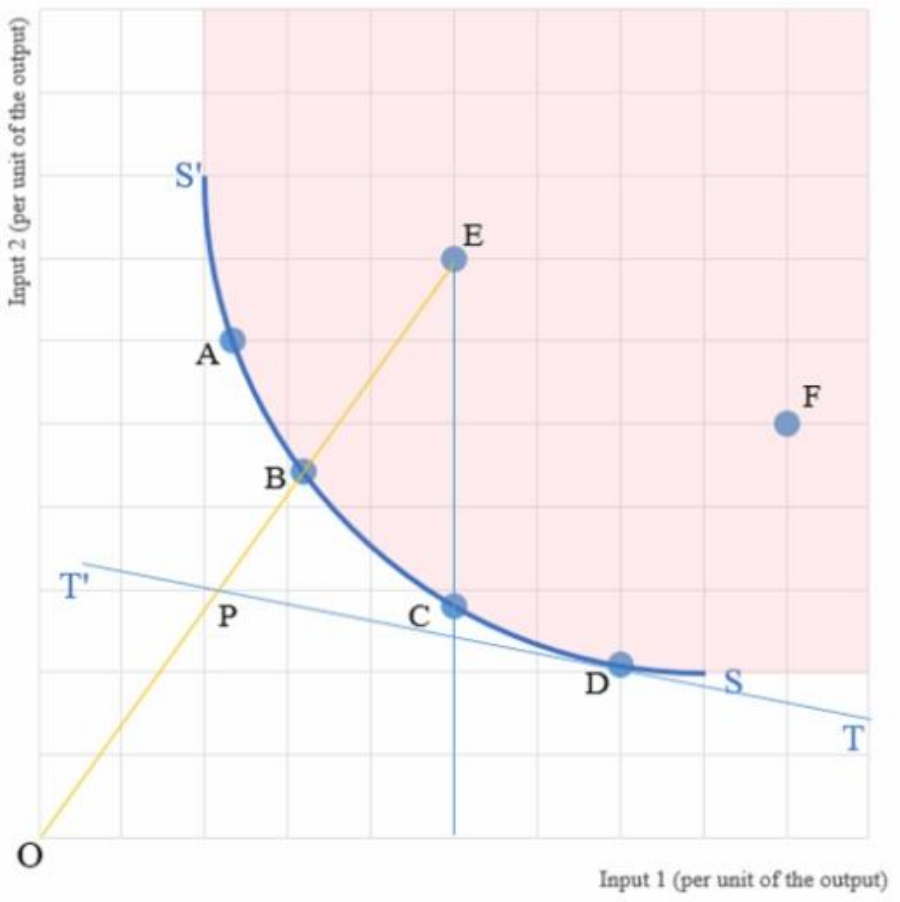

Figure 2: The price/overall/relative/allocative/economic efficiency of E.

From the literature of economics, the ratio of $OB/OE$ is called *the technical efficiency of E*, the ratio of $OP/OB$ is called *the price efficiency of B* or *the allocative efficiency* of E, and the ratio of $OP/OE$ is called the ratio of *the overall efficiency of E* or *the economic efficiency* of E. All of these ratios are less than equal to 1.

As can be seen, the overall efficiency of E has exactly the same meaning as *the price efficiency of E*, which can also be measured by multiplying the technical efficiency of E and the allocative efficiency of E, that is, $(OB/OE) \times (OP/OB) = OP/OE$. Therefore, at least two of phrases, 'price efficiency', 'overall efficiency' or 'economic efficiency', are redundant.

The concept of allocative efficiency also does not provide a ranking tool similar to the concept of technical efficiency and just describes *the non-technical efficiency*; hence, there is no reason to use a new phrase. Indeed, *the price inefficiency of E* can be decomposed by *the technical inefficiency of E* and *the non-technical inefficiency of E*, thus we should avoid using extra and superfluous phrases.

From these definitions, the price (overall/economic) efficiency of D is 1 and the price efficiency of the other firms is always less than equal to 1, thus, the meaning of the price (overall/economic) efficiency can also be interpreted as 'the relative efficiency'. The equation $w_1x_1 + w_2x_2$ has the same value for every point $(x_1, x_2)$ on the line $TT'$, and is equal to $w_1 x_{1D} + w_2 x_{2D}$, where $w_1$ and $w_2$ are the corresponded prices/weights/worth of the first input factor and the second input factor, respectively.

As Figure 3 illustrates, the price (overall/economic) efficiency of E, $OP/OE$, is equal with the ratio of $(w_1 x_{1P} + w_2 x_{2P})/(w_1 x_{1E} + w_2 x_{2E})$, which is equal to $(w_1 x_{1D} + w_2 x_{2D})/(w_1 x_{1E} + w_2 x_{2E})$. Indeed, the lines $TT'$ and $RR'$ are parallel and the ratio of $OP/OE$ is equal to the ratio of $OP'/OE'$ in the triangle $OEE'$. Hence, the linear combination of the input factors of every point of the PPS is compared with the linear combination of D's input factors. Thus, the provided score is relatively meaningful, and the price (overall/economic) efficiency of E can be introduced as the relative efficiency of E as well. This outcome precisely defines the concept of doing the job well and can be expressed in one word 'efficiency', which is the condition or fact of producing the results that we would want without waste.

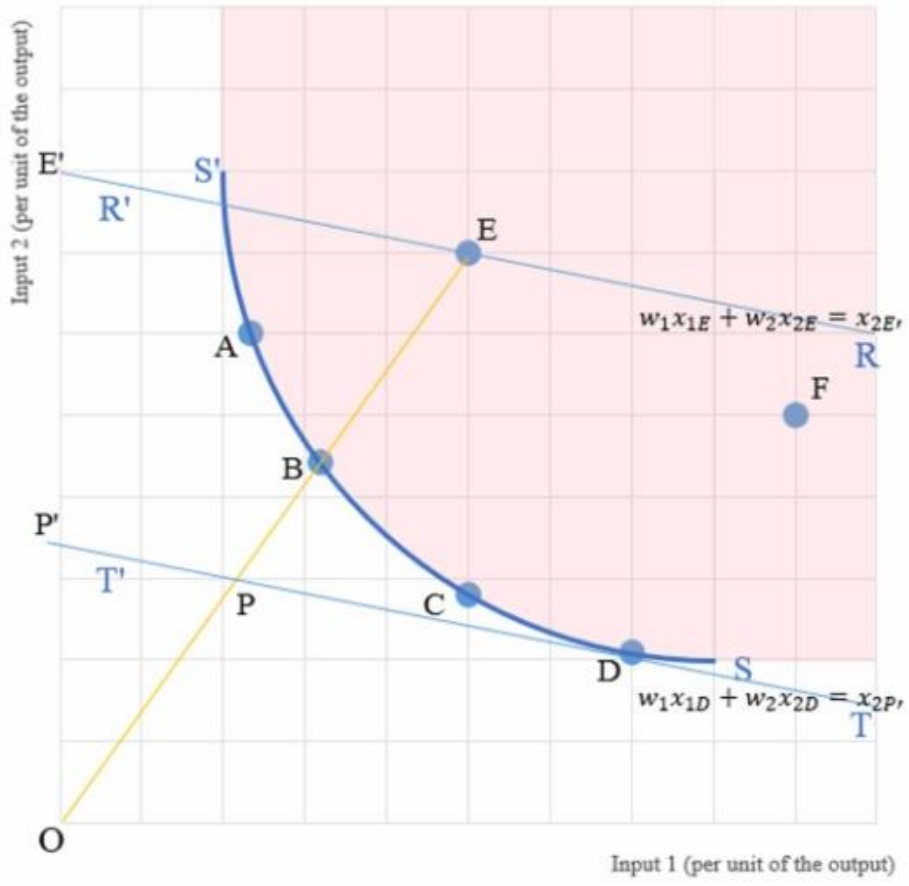

Figure 3: The concept of doing the job well.

Note that, in the literature of operations research, 'relative efficiency' is usually considered as 'technical efficiency', which is incorrect. In definition of technical efficiency, no suitable discrimination between the points on the production function is introduced. The provided ratio for the technical efficiency, such as, $OB/OE$, is a fake relative score, and does not yield a valid comparison between E and other points in the PPS. In fact, the pros and cons of the technical efficiency are the same as that of doing the job right, and the provided score for the concept of doing the job right is not relatively meaningful.

Since, the concept of doing the job well depends on the weights/prices/worth of the factors, and requires the concepts of the wholly dominant, the convexity and the radiate approaches to discriminate the firms linearly; the efficiency also depends on the weights/prices/worth of the factors, and at least requires the concepts of the wholly dominant, the convexity and the radiate approaches to discriminate the firms linearly.

In the literature of economics and operations research, the wholly dominant approach is called *Free Disposal Hull* (FDH) technology, the combination of the wholly dominant and the convexity approaches is called *Variable Returns to Scale* (VRS) technology, and the combination of the wholly dominant, the convexity and the radiate approaches is called *Constant Returns to Scale* (CRS) technology, as Table 1 illustrates. The technical efficiency (doing the job right) depends on the FDH, VRS or CRS technologies and does not provide discrimination between firms, but the efficiency (doing the job well) depends on the relationships between the factors and at least requires *the CRS technical efficiency* to discriminate firms linearly.

When a firm is not efficient, it is *inefficient*. If one desires to decompose the inefficiency of a firm, inefficiency can be decomposed by the CRS technical inefficiency and the non-CRS-technical inefficiency. The CRS technical inefficiency can also be decomposed by VRS technical inefficiency and non-VRS technical inefficiency, and so on. This topic is discussed in the upcoming sections.

In short, the meaning of CRS technical efficiency should not be misinterpreted as efficiency, similar to the concept of doing the job right which should not be misinterpreted with the concept of doing the job well. If using $200 at most yields $200, and $220 yields $700, the point (220, 700) is more efficient than the point (200, 200), and this is our suitable choice, regardless of whether we are applying VRS, FDH or any other approaches to define the production function. In other words, when our purpose is to rank a set of homogenous firms, at least CRS technical efficiency should be measured. Of course, after finding the best firm and measure the concept of partially dominant, the exact returns to scale is required to estimate the production function. *From here forward, instead of the phrase 'CRS technical efficiency' the phrase 'technical efficiency' will be exclusively used.*

## 2.3. The productivity measurement

The concept of doing the job well is also introduced as '*productivity*' in the literature of economics. There is no problem if one desires to call 'doing the job well' as productivity and one can use the pair 'technical efficiency and productivity' instead of the pair 'technical efficiency and efficiency'; nonetheless, after measuring the concept of doing the job well, there is still a need to measure whether the outcomes satisfy the goals of firms. Indeed, the concept of 'doing the well job' is different from the concept of 'doing the job well', and requires another meaningful name.

According to Cambridge English dictionary, "the ability to be successful and produce the intended results" is called '*effectiveness*'. Therefore, 'effectiveness' can be used for the concept of 'doing the well job' and 'productivity' can be used for the concept of 'doing the useful job', which is a combination of both efficiency and effectiveness. The word 'productivity' means "the rate at which a person, company or country does useful work", according to Cambridge English dictionary. Therefore, it is suggested that the commonly utilized phrases and concepts be reestablished according to the following table.

Table 1: Reestablishing the concepts/phrases.

| The concept/approach | Reestablishing |
|---|---|
| Doing the job right | Technical efficiency |
| Doing the job well | Efficiency |
| Doing the right job | Technical effectiveness |
| Doing the well job | Effectiveness |
| Doing the useful job | Productivity |
| The wholly dominant approach | Free disposal hull |
| The wholly dominant and convexity approaches | Variable returns to scale |
| The wholly dominant, convexity and radiate approaches | Constant returns to scale |
| The wholly dominant, convexity and inner radiate approaches | Decreasing returns to scale |
| The wholly dominant, convexity and outer radiate approaches | Increasing returns to scale |

To clearly explain the above, let's suppose that a set of 9 homogenous banks, labeled A-I, are selected. Assume that the aim of discrimination is (1) to find the banks which have used a smaller number of tellers to service a greater number of customers, and (2) to find the banks which have at least serviced $y_1$ number of customers in the period of evaluation. Suppose that the production function is available and the location of each bank in the Cartesian coordinate plane is depicted in Figure 4.

The blue curve represents the production function, the horizontal axis illustrates the number of tellers and the vertical axis displays the number of customers. The banks A-G lie on the production function, and are VRS-technically efficient, and the banks H and I are VRS-technically inefficient.

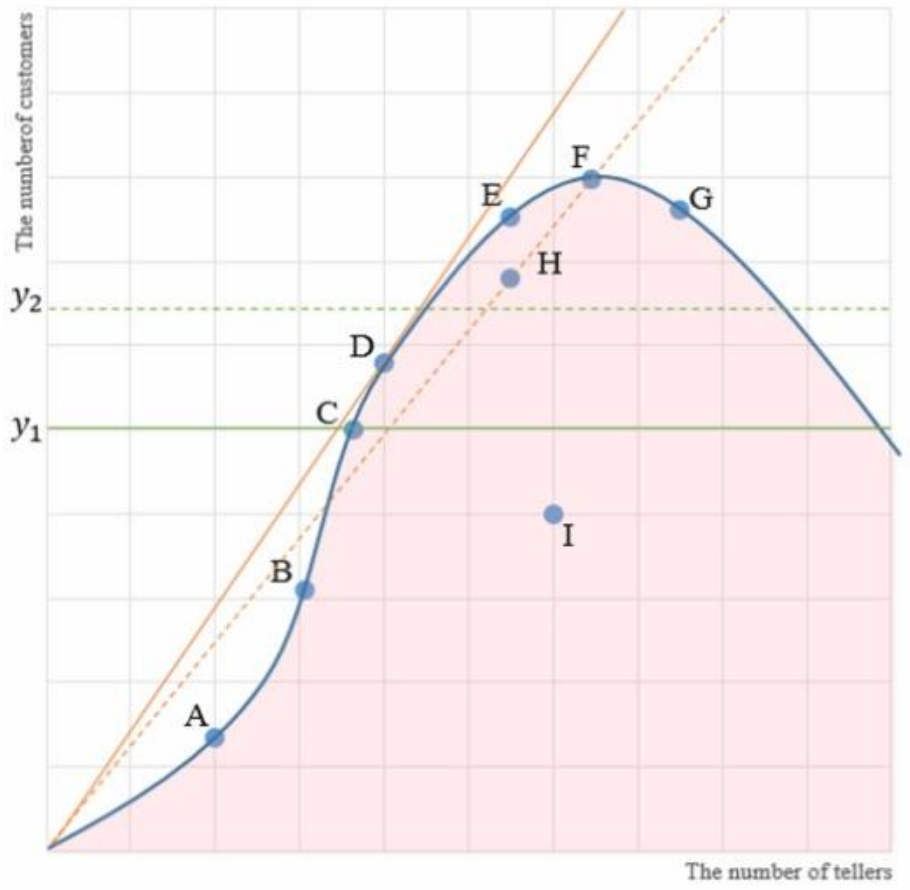

Figure 4: The production function.

By considering the ratio of the number of customers to the number of tellers, D is the most efficient bank followed by E, H and F, respectively. H is not technically efficient, but, for instance, H is more efficient than A which is VRS-technically efficient.

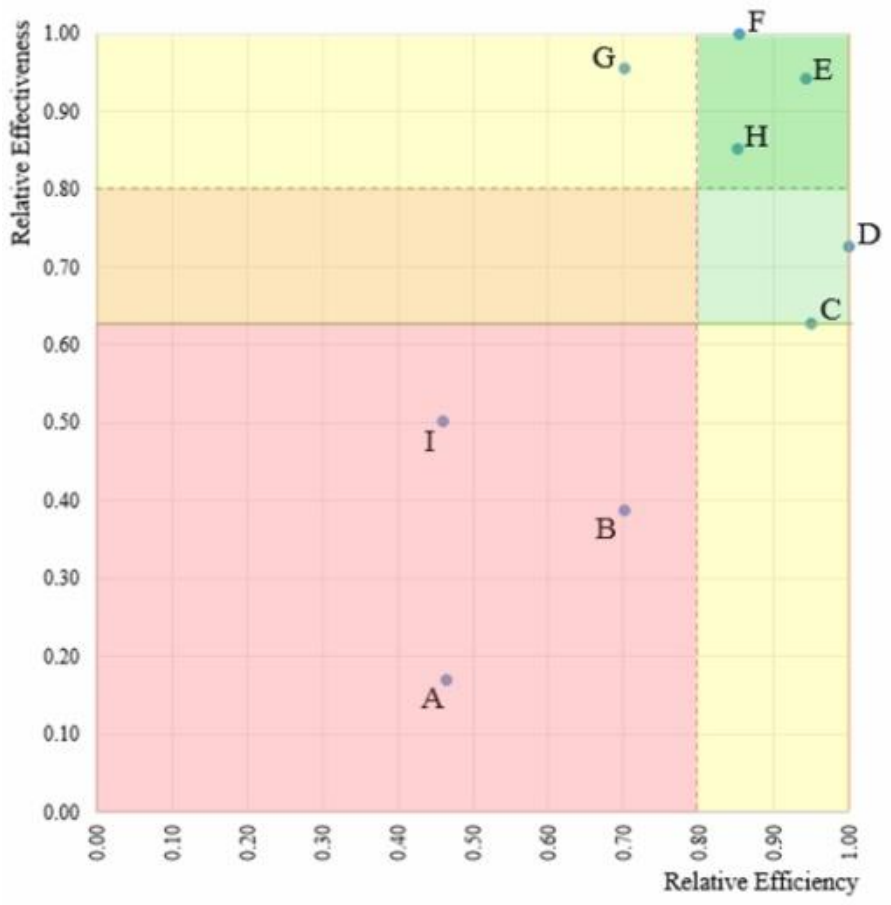

Figure 5: The productivity measurement.

The banks C, D, E, F, G and H have at least serviced $y_1$ number of customers, and are *effective*. F is the most effective bank followed by G, E and H. Therefore, the banks C, D, E, F and H are the most productive banks which are most efficient and most effective at the same time in comparison with other banks.

The above results can also be seen in Figure 5, where the horizontal ax-is displays the relative efficiency of each bank and the vertical axis represents the relative effectiveness of each bank. As can be seen, D has the relative score equal to 1 and F has the relative effectiveness equal to 1. The red area illustrates the non-productive banks, which are the banks A, B and I, and have relative effectiveness scores less than that of C and relative efficiency scores less than 0.8.

G is the most effective bank after F, but has the relative efficiency score less than 0.8. The (dark and light) green area represents the productive banks which have relative efficiency scores more than 0.8 and relative effectiveness scores more than that of C.

If the goal of evaluation, which is at least servicing $y_1$ number of customers in period of evaluation, is changed to at least servicing $y_2$ number of customers, as Figures 4 and 5 illustrate, even the most efficient banks D and C are not called productive due to the lack of their effectiveness. In this case, the banks E, F and H are the most productive banks among the banks A-I, as the dark green area displays in Figure 5.

The concept of effectiveness is always required in real-life applications, for instance, no bank desires to decrease *consumer satisfaction* and no firm works without a plan or some requirements. It is also possible that each bank has different goals for effectiveness index; for instance, a bank may only prefer the production function values according to its set of input factors and in this case, if a different goal for each bank does not affect the homogeneity of the banks, the technical effectiveness index can also be calculated.

In short, there are several concepts which provide the most important indexes to discriminate the performance of homogenous firms. The concept of doing the job right, which is called technical efficiency, is not enough to discriminate between firms. The technical efficiency is usually interpreted as efficiency in the literature, and if one would like to use such a term for the concept of doing the job right, one should be aware that such efficiency is neither enough to discriminate between firms, nor should the corresponded index be used to rank and benchmark the firms. In order to rank and benchmark a set of homogenous firms, the concept of doing the job well is required, and this concept is called efficiency. The efficiency is usually interpreted as productivity in the literature of economics, and similarly if one would like to use the word productivity for the concept of doing the job well, one should know that there is still a need to measure the effectiveness of firms. Therefore, Table 1 is provided to reintroduce the concepts of efficiency, effectiveness and productivity with industry-wide accuracy and understanding.

## 3.0. The epsilon KAM

Suppose that there are $n$ firms, labeled $F_i$ $(i = 1,2,\dots,n)$, and each firm has $m$ input factors with the values $x_{ij}$ $(j = 1,2,\dots,m)$ and $p$ output factors with the values $y_{ik}$ $(k = 1,2,\dots,p)$. Assume that the weights/prices or the approximation of the relationships between input and output factors are $W_j^-$ and $W_k^+$, for $j = 1,2,\dots,m$ and $k = 1,2,\dots,p$, respectively. Suppose that, $V_j^-$ and $V_k^+$ are defined as Equation 1, for $j = 1,2,\dots,m$ and $k = 1,2,\dots,p$:

$$V_j^- = \frac{W_j^-}{\sum_{j=1}^{m} W_j^- x_{lj}} \quad \& \quad V_k^+ = \frac{W_k^+}{\sum_{k=1}^{p} W_k^+ y_{lk}}. \tag{1}$$

Assume that the epsilon vector, with bolded notation $\boldsymbol{\epsilon}$, is given by $\boldsymbol{\epsilon} = (\varepsilon_1^-, \varepsilon_2^-, \ldots, \varepsilon_m^-, \varepsilon_1^+, \varepsilon_2^+, \ldots, \varepsilon_p^+)$, to introduce an epsilon neighborhood of firm $F_l$, $(l = 1,2,\ldots,n)$. The components of epsilon vector are introduced by Equation 2, for $j = 1,2,\ldots,m$ and $k = 1,2,\ldots,p$.

$$\varepsilon_j^- = \varepsilon / W_j^- \quad \text{and} \quad \varepsilon_k^+ = \varepsilon / W_k^+. \tag{2}$$

The value of epsilon has the same meaning for each factor, when Equation 2 is considered.

The components of epsilon vector can be defined by Equation 3, where $\varepsilon \in [0, +\infty)$, for $j = 1,2,\ldots,m$, and $k = 1,2,\ldots,p$,

$$\varepsilon_j^- = \varepsilon \times x_{lj} \quad \text{and} \quad \varepsilon_k^+ = \varepsilon \times y_{lk}. \tag{3}$$

There are a lot of ways to introduce the components of epsilon vector, according to the aim of discrimination. For instance, $\varepsilon_j^-$ and $\varepsilon_k^+$ can be introduced as Equation 4, for $j = 1,2,\ldots,m$, and $k = 1,2,\ldots,p$. In this case, the components of epsilon vector are commensurate with the corresponded input and output factors, but are not changed from one firm to another.

$$\varepsilon_j^- = \varepsilon \times \operatorname*{ave}_{1 \le i \le n} x_{ij} \quad \text{and} \quad \varepsilon_k^+ = \varepsilon \times \operatorname*{ave}_{1 \le i \le n} y_{ik}. \tag{4}$$

It is also possible to introduce one (or more) of the components of epsilon vector as 0, regarding the purpose of discrimination. For example, the components of epsilon can be introduced by $\varepsilon_j^- = \varepsilon \times x_{lj}$ and $\varepsilon_k^+ = 0$, for $j = 1,2,\ldots,m$, and $k = 1,2,\ldots,p$, which let's consider the errors in input factors only.

Equation 5 illustrates the $\boldsymbol{\epsilon}$-KAM, when the performance of $F_l$ is measured $(l = 1,2,\ldots,n)$.

$$\min \frac{\sum_{j=1}^{m} V_j^- (x_{lj} + \varepsilon_j^- - s_j^-)}{\sum_{k=1}^{p} V_k^+ (y_{lk} - \varepsilon_k^+ + s_k^+)}, \tag{5}$$

Subject to

$\sum_{i=1}^{n} \lambda_i x_{ij} + s_j^- = x_{lj} + \varepsilon_j^-$, for $j = 1,2,\ldots,m$,

$\sum_{i=1}^{n} \lambda_i y_{ik} - s_k^+ = y_{lk} - \varepsilon_k^+$, for $k = 1,2,\ldots,p$,

$\lambda_i \ge 0$, for $i = 1,2,\ldots,n$,

$s_j^- \ge 0$, for $j = 1,2,\ldots,m$,

$s_k^+ \ge 0$, for $k = 1,2,\ldots,p$.

The $\boldsymbol{\epsilon}$-KAM can linearly be solved by Equation (6).

$$\min\left[\sum_{j=1}^{m} V_j^- (t x_{lj} + t\varepsilon_j^- - s_j^-)\right], \tag{6}$$

Subject to

$\left[\sum_{k=1}^{p} V_k^+ (t y_{lk} - t\varepsilon_k^+ + s_k^+)\right] = 1$,

$\sum_{i=1}^{n} \lambda_i x_{ij} + s_j^- = t x_{lj} + t\varepsilon_j^-$, for $j = 1,2,\ldots,m$,

$\sum_{i=1}^{n} \lambda_i y_{ik} - s_k^+ = t y_{lk} - t\varepsilon_k^+$, for $k = 1,2,\ldots,p$,

$\lambda_i \ge 0$, for $i = 1,2,\ldots,n$,

$s_j^- \ge 0$, for $j = 1,2,\ldots,m$,

$s_k^+ \ge 0$, for $k = 1,2,\ldots,p$,

$t > 0$.

The score of KAM represents that the efficiency score of firm $l$, that is, $\sum_{k=1}^{p} W_k^{+} y_{lk} / \sum_{j=1}^{m} W_j^{-} x_{lj}$, is compared with the efficiency score of a point on the estimated production function , that is, $\sum_{k=1}^{p} W_k^{+}(y_{lk} - \varepsilon_k^{+} + s_k^{+}) / \sum_{j=1}^{m} W_j^{-}(x_{lj} + \varepsilon_j^{-} - s_j^{-})$, such that, the ratio of $\sum_{k=1}^{p} W_k^{+} y_{lk} / \sum_{j=1}^{m} W_j^{-} x_{lj}$ to $\sum_{k=1}^{p} W_k^{+}(y_{lk} - \varepsilon_k^{+} + s_k^{+}) / \sum_{j=1}^{m} W_j^{-}(x_{lj} + \varepsilon_j^{-} - s_j^{-})$ becomes minimum. Since the efficiency of firm $l$ is a constant value, KAM finds a point on the estimated production function which has an equal or greater efficiency score in comparison with that of firm $l$, regarding the value of epsilon and introduced $W_j^{-}$ and $W_k^{+}$.

The target for firm $l$ ($l = 1,2,\dots,n$), which lies on the estimated production function, and has a greater (or equal) ratio of the linear combination of output factors to the linear combination of input factors, regarding the value of epsilon and introduced $W_j^{-}$ and $W_k^{+}$, is given by:

$$
\begin{aligned}
&x_{lj}^{*} = x_{lj} + \varepsilon_j^{-} - s_j^{-*}/t^{*}, \text{ for } j = 1,2,\dots,m, \\
&y_{lk}^{*} = y_{lk} - \varepsilon_k^{+} + s_k^{+*}/t^{*}, \text{ for } k = 1,2,\dots,p.
\end{aligned} \tag{7}
$$

The dual linear programming is also given by Equation 8.

$$
\begin{aligned}
&\max \tau, \\
&\tau\left(\textstyle\sum_{k=1}^{p} V_k^{+}(y_{lk} - \varepsilon_k^{+})\right) + \textstyle\sum_{j=1}^{m} w_j^{-}\left(x_{lj} + \varepsilon_j^{-}\right) - \\
&\textstyle\sum_{k=1}^{p} w_k^{+}(y_{lk} - \varepsilon_k^{+}) = \textstyle\sum_{j=1}^{m} V_j^{-}(x_{lj} + \varepsilon_j^{-}), \\
&\textstyle\sum_{k=1}^{p} w_k^{+} y_{ik} - \textstyle\sum_{j=1}^{m} w_j^{-} x_{ij} \le 0, \text{ for } i = 1,2,\dots,n, \\
&w_j^{-} \ge V_j^{-}, \text{ for } j = 1,2,\dots,m, \\
&w_k^{+} \ge \tau V_k^{+}, \text{ for } k = 1,2,\dots,p.
\end{aligned} \tag{8}
$$

Suppose that the component of epsilon vector is introduced by Equation 2. When $\varepsilon = 0$, the 0-KAM measures the technical efficiency of firms, and divides the firms into two categories similar to DEA models. If the score of 0-KAM is less than 1 for a firm, that firm is technically inefficient, and if the score of 0-KAM is 1, the firm is technically efficient. However, the scores of 0-KAM (similar to DEA models) should neither be used to rank firms, nor the proposed targets can be used to benchmark firms.

The $\varepsilon$-KAM is SBM (Tone 2001), where $\varepsilon = 0$, $W_j^{-} = 1/x_{lj}$ and $W_k^{+} = 1/y_{lk}$. When $\varepsilon > 0$, and $W_j^{-} = 1/x_{lj}$ and $W_k^{+} = 1/y_{lk}$, we express '*the $\varepsilon$-KAM with SBM approach*'. Note that, the SBM approach does not measure efficiency, but can be used to introduce the technically efficient firms. When $\varepsilon > 0$, and $W_j^{-} = 1/\min\{x_{lj}: x_{lj} \ne 0\}$ and $W_k^{+} = 1/\min\{y_{lk}: y_{lk} \ne 0\}$, we state '*the $\varepsilon$-KAM with minimum approach*'.

Suppose that $W_j^{-}$ and $W_k^{+}$ are given as the available costs of input and output factors, where $j = 1,2,\dots,m$ and $k = 1,2,\dots,p$. We can decrease the discrimination power of $\varepsilon$-KAM, as Equations 9, 10 represent, to illustrate the lack of CF and RF (or PF) measurements, respectively.

$$
\begin{aligned}
&\min \textstyle\sum_{j=1}^{m} V_j^{-}(x_{lj} + \varepsilon_j^{-} - s_j^{-}), \\
&\text{Subject to} \\
&\textstyle\sum_{i=1}^{n} \lambda_i x_{ij} + s_j^{-} = x_{lj} + \varepsilon_j^{-}, \text{ for } j = 1,2,\dots,m, \\
&\textstyle\sum_{i=1}^{n} \lambda_i y_{ik} - s_k^{+} = y_{lk}, \text{ for } k = 1,2,\dots,p, \\
&\lambda_i \ge 0, \text{ for } i = 1,2,\dots,n, \\
&s_j^{-} \ge 0, \text{ for } j = 1,2,\dots,m, \\
&s_k^{+} \ge 0, \text{ for } k = 1,2,\dots,p.
\end{aligned} \tag{9}
$$

$$\max \sum_{k=1}^{p} V_k^+ (y_{lk} - \varepsilon_k^+ + s_k^+), \tag{10}$$

Subject to

$\sum_{i=1}^{n} \lambda_i x_{ij} + s_j^- = x_{lj}$, for $j = 1,2,\dots,m$,

$\sum_{i=1}^{n} \lambda_i y_{ik} - s_k^+ = y_{lk} - \varepsilon_k^+$, for $k = 1,2,\dots,p$,

$\lambda_i \geq 0$, for $i = 1,2,\dots,n$,

$s_j^- \geq 0$, for $j = 1,2,\dots,m$,

$s_k^+ \geq 0$, for $k = 1,2,\dots,p$.

In other words, none of Equations 9 and 10 provide a fair measure to discrimination the efficiency of firms. A suitable model should at least satisfy the following Types 1-6. For instance, Equation 9 does not satisfy Type 4, and Equation 10 does not satisfy Type 3.

A firm, $F_l$, partially dominates another firm, $F_{l'}$, if and only if, the value of Equation 11 for $F_l$ is greater than that of $F_{l'}$, where $l$ and $l'$ belongs to $\{1,2,\dots,n\}$.

$$\frac{\sum_{k=1}^{p} W_k^+ y_k}{\sum_{j=1}^{m} W_j^- x_j} = \frac{W_k^+ y_1 + W_k^+ y_2 + \cdots + W_k^+ y_p}{W_j^- x_1 + W_j^- x_2 + \cdots + W_j^- x_m}. \tag{11}$$

The six introduced types which increase the value of Equation 11 are expressed by:

**Type 1**: Decreasing the value of denominator in Equation 11 when the value of numerator is fixed.

**Type 2**: Increasing the value of numerator in Equation 11 when the value of denominator is fixed.

**Type 3**: If the rate of increasing the value of numerator in Equation 11 is greater than the rate of increasing the value of denominator.

**Type 4**: If the rate of decreasing the value of numerator in Equation 11 is greater than the rate of decreasing the value of denominator.

**Type 5**: The value of denominator in Equation 11 is decreased (increased) by: (1) decreasing (increasing) the value of one or more of input factors, or (2) increasing (decreasing) a small value of one or more of input factors and decreasing (increasing) a large value of one or more of other input factors.

**Type 6**: The value of numerator in Equation 11 is increased (decreased) by: (1) increasing (decreasing) the value of one or more of output factors, or (2) decreasing (increasing) a small value of one or more of output factors and increasing (decreasing) a large value of one or more of other output factors.

In short, Equation 5 shows that KAM compares the efficiency of a firm with the efficiency of the points on the estimated production function, and finds the best target for the firm, regarding the value of epsilon and introduced weights. The discrimination power of KAM is greater than CF, RF and PF models, and KAM provides a complete measure to assess all the inefficiencies of the firms. In the next sections, KAM is improved, and the optimum of epsilon is also measured.

## 4.0. KAM and uncontrollable factors

A researcher wants to compare eight homogenous international Iranian airports, labeled A-H, according to seven factors, the area of airport (Hectare), the area of apron (Square meter), the area of terminal (Square meter), the area of runway (Square meter), the number of operating flights, the number of passengers' movements, the amount of air cargo (Metric ton), as Table 2 represents. Which airport does the job well?

Table 2. The data of 8 homogenous airports.

| N | Airport | Area | Apron | Terminal | Runway | Flights | Passengers | Cargo |
|---|---|---|---|---|---|---|---|---|
| 1 | A | 1,200 | 304,182 | 45,600 | 353,610 | 30,707 | 4,030,859 | 74,184 |
| 2 | B | 503 | 213,729 | 38,778 | 348,120 | 46,875 | 4,783,120 | 19,050 |
| 3 | C | 800 | 41,003 | 11,800 | 269,955 | 15,608 | 1,039,967 | 1,587 |
| 4 | D | 1,041 | 112,464 | 21,050 | 395,730 | 39,871 | 1,744,524 | 4,919 |
| 5 | E | 1,002 | 30,000 | 8,000 | 192,330 | 4,887 | 427,974 | 1,574 |
| 6 | F | 478 | 63,000 | 23,000 | 389,115 | 41,088 | 2,165,572 | 5,414 |
| 7 | G | 481 | 47,210 | 9,300 | 268,995 | 19,010 | 971,313 | 3,826 |
| 8 | H | 1,346 | 503,274 | 76,370 | 421,305 | 129,153 | 11,709,741 | 39,556 |

The first four factors, the areas of airport, apron, terminal and runway, are input factors, because they illustrate the infrastructure of airports and lesser values of these factors have worth. The last three factors, the numbers of operating flights and passengers' movements and the amounts of air cargo, are output factors, because they represent the business and production of airports and greater values of these factors have worth. An airport which uses lesser amounts of input factors and greater amounts of output factors has done the job well in comparison with another airport.

The runway is an input factor which has the standard area, according to the documents of International Civil Aviation Organization (ICAO). The area of runway should not be less than the standard value and depends on the aircrafts which want to land and take off from that runway. Because of the safety of passengers, decreasing the length/area of a runway is not suggested. This kind of factor is called *non-controllable* or *uncontrollable* factor (Banker and Morey, 1986; Charnes *et al*., 1987; Cooper *et al*., 2007). In other words, an uncontrollable factor may not be controlled by managers, although, it may affect the performance of firms.

For instance, suppose that $j^{\text{th}}$ input factor ($j = 1,2,\dots,m$) is an uncontrollable factor and the efficiency of $F_l$ is measured ($l = 1,2,\dots,n$). If this uncontrollable factor should not be decreased and increased, the corresponded linear combination of this $j^{\text{th}}$ input factor of firms in Equation 4.7, (that is, $\sum_{i=1}^{n}\lambda_i x_{ij}$) should be equal to the corresponded value of $j^{\text{th}}$ factor of firm $l$ (that is, $x_{lj}$), that is, $\sum_{i=1}^{n}\lambda_i x_{ij} = x_{lj}$, (where $l = 1,2,\dots,n$). Nonetheless, it is valuable to examine an epsilon neighborhood of an uncontrollable factor, to measure the effect of such restriction on other factors. For such an aim, the constraint $s_j^- \leq \varepsilon_j^-$ (or $s_j^- \leq t\varepsilon_j^-$) can be added to the constraints in Equation 5 (Equation 6).

The term $\varepsilon_j^-$ is a value which is added to $j^{\text{th}}$ input factor of $F_l$, to introduce a neighborhood of this factor, thus the linear combination of this uncontrollable factor of firms is at least equal to $x_{lj}$, that is, $\sum_{i=1}^{n}\lambda_i x_{ij} = x_{lj} + \varepsilon_j^- - s_j^- \geq x_{lj}$. Therefore, the optimal value of $j^{\text{th}}$ input factor is not decreased, but it may slightly be increased according to value of $\varepsilon_j^-$.

For the case that the optimal value of $j^{\text{th}}$ input factor should not also be increased, the value of epsilon can be considered very small, such that, the value of $\varepsilon_j^-$ is quite negligible. For instance, if $j^{\text{th}}$ input factor is measured with three decimal digits, the negligible error can be less than 0.005. Similar discussion can be illustrated for an uncontrollable output factor as well.

Now, suppose that $J_u$ is a subset of input factor indexes, $\{1,2,\dots,m\}$, corresponded to the uncontrollable input factors, and $K_u$ is a subset of output factor indexes, $\{1,2,\dots,p\}$, corresponded to the uncontrollable output factors. The $\varepsilon$-KAM is given by Equation 12 (See also Khezrimotlagh *et al* 2012a,b, 2013 and 2014).

$$\begin{aligned}
&\min\left[\sum_{j=1}^{m} V_j^-(tx_{lj} + t\varepsilon_j^- - s_j^-)\right], \qquad (12)\\
&\text{Subject to}\\
&\left[\sum_{k=1}^{p} V_k^+(ty_{lk} - t\varepsilon_k^+ + s_k^+)\right] = 1,\\
&\sum_{i=1}^{n} \lambda_i x_{ij} + s_j^- = tx_{lj} + t\varepsilon_j^-, \text{ for } j = 1,2,\dots,m,\\
&\sum_{i=1}^{n} \lambda_i y_{ik} - s_k^+ = ty_{lk} - t\varepsilon_k^+, \text{ for } k = 1,2,\dots,p,\\
&s_j^- \le t\varepsilon_j^-, \text{ for } j \in J_u,\\
&s_k^+ \ge t\varepsilon_k^+, \text{ for } k \in K_u,\\
&\lambda_i \ge 0, \text{ for } i = 1,2,\dots,n,\\
&s_j^- \ge 0, \text{ for } j = 1,2,\dots,m,\\
&s_k^+ \ge 0, \text{ for } k = 1,2,\dots,p,\\
&t > 0.
\end{aligned}$$

For an example of uncontrollable factor, suppose that fourth input factor (runway) is uncontrollable in the above airport example. Assume that $W_j^- = 1/x_{lj}$, for $j = 1,2,3,4$, $W_k^+ = 1/y_{lk}$, for $k = 1,2,3$, where $\varepsilon$ is 0, 0.0001, 0.01 and 0.1, and the component of epsilon vectors are introduced by Equation 3. Since the weights are introduced by the inverse of data like SBM, the approach is the SBM approach and the results are only used to express the methodology.

When $\varepsilon = 0$, in Equation 12, the corresponded constraints for fourth input factor of firm $F_l$, (that is, $s_4^- \le t\varepsilon x_{l4}$, and $\sum_{i=1}^{8} \lambda_i x_{i4} + s_4^- = t(1+\varepsilon)x_{l4}$) are equal to zero, that is, $s_4^- = t\varepsilon x_{l4} = 0$, and $\sum_{i=1}^{8} \lambda_i x_{i4} = tx_{l4}$. As Table 3 illustrates, the 0-KAM by SBM approach represents that all airports A-H are technically efficient. In other words, this simple restriction on fourth factor of airports yields that airports C and E become technically efficient.

Table 3: The $\varepsilon$-KAM scores where $s_4^- = t\varepsilon x_{l4}$.

| Epsilon | A | B | C | D | E | F | G | H |
|---|---|---|---|---|---|---|---|---|
| 0 | 1.00000 | 1.00000 | 1.00000 | 1.00000 | 1.00000 | 1.00000 | 1.00000 | 1.00000 |
| 0.0001 | 0.99983 | 0.99997 | 0.99924 | 0.98588 | 0.99258 | 0.99993 | 0.99990 | 0.99999 |
| 0.01 | 0.98307 | 0.99747 | 0.92921 | 0.63590 | 0.54534 | 0.99274 | 0.98991 | 0.99857 |
| 0.1 | 0.85672 | 0.97644 | 0.55163 | 0.50990 | 0.36364 | 0.93414 | 0.91129 | 0.98565 |

As the value of epsilon is increased, (for instance, when $\varepsilon$ is 0.0001, 0.01 and 0.1), $\varepsilon$-KAM discriminates the airports according to the introduced value of epsilon and the SBM approach.

If the uncontrollable factor cannot be decreased and increased, the constraint $s_4^- \le t\varepsilon x_{l4}$ can be replaced by $s_4^- = t\varepsilon x_{l4}$. The results of $\varepsilon$-KAM for this assumption are illustrated in Table 3. In this case, the value of fourth factor for airport number $l$ is not changed. Indeed, the targets for firm $l$ $(l = 1,2,\dots,n)$ can be measured from Equation 7. Since $s_4^- = t\varepsilon_j^- = t\varepsilon x_{l4}$, for instance, when $\varepsilon = 0.0001$, Table 4 represents that the suggested target for fourth input factor of airport number $l$ is not changed, $(l =$

1,2, … ,8).

Table 4: The $10^{-4}$-KAM targets where $s_4^- = t\varepsilon x_{l4}$.

| Airport | Area | Apron | Terminal | Runway | Flights | Passengers | Cargo |
|---|---|---|---|---|---|---|---|
| A | 1,187.28 | 325,581.90 | 48,948.45 | 353,610.00 | 44,770.19 | 5,080,228.95 | 66,765.60 |
| B | 553.30 | 230,417.98 | 40,786.56 | 348,120.00 | 51,816.17 | 5,187,091.37 | 20,175.81 |
| C | 418.46 | 45,103.30 | 12,980.00 | 269,955.00 | 23,914.76 | 1,258,627.24 | 3,681.02 |
| D | 778.22 | 123,710.40 | 23,155.00 | 395,730.00 | 35,883.90 | 2,441,589.64 | 17,250.27 |
| E | 297.08 | 33,000.00 | 8,800.00 | 192,330.00 | 16,583.00 | 865,955.98 | 2,825.11 |
| F | 497.53 | 69,300.00 | 23,630.14 | 389,115.00 | 40,918.80 | 2,218,202.93 | 7,181.15 |
| G | 486.49 | 51,931.00 | 10,230.00 | 268,995.00 | 19,494.13 | 1,047,313.36 | 5,212.20 |
| H | 1,357.68 | 483,620.85 | 73,294.86 | 421,305.00 | 116,237.70 | 10,746,024.81 | 46,368.88 |

The suggested targets in Table 4 (or the targets by Equation 5) lie on the frontier of feasible area and have better performance in comparison with the real data in Table 2, regarding the assumptions of discrimination. In other words, 0.0001-KAM not only discriminates the airports according to the introduced errors, but it also benchmarks all technically efficient airports, and suggests how they can regulate their factors according to Types 1-6, to improve their performances, (regarding the introduced assumptions). For instance, 0.0001-KAM says that A should decrease the area of airport and increase the numbers of flights and passengers, even if the area of apron and terminal are increased, or the amount of cargo decreased, according to the value of epsilon and the SBM approach.

Now, suppose that $s_4^- \leq t\varepsilon x_{l4}$ is only added, the corresponded results to Tables 3 are displayed in Table 5. The constraint, $s_4^- \leq t\varepsilon x_{l4}$, does not let the value of fourth factor decrease, but it may be increased to find a better situation on the production frontier, according to the SBM approach.

Table 5: The $\varepsilon$-KAM scores where $s_4^- \leq t\varepsilon x_{l4}$.

| Epsilon | A | B | C | D | E | F | G | H |
|---|---|---|---|---|---|---|---|---|
| 0 | 1.00000 | 1.00000 | 1.00000 | 1.00000 | 1.00000 | 1.00000 | 1.00000 | 1.00000 |
| 0.0001 | 0.99980 | 0.99997 | 0.99924 | 0.98468 | 0.99258 | 0.99993 | 0.99990 | 0.99997 |
| 0.01 | 0.98056 | 0.99747 | 0.92921 | 0.63590 | 0.54534 | 0.99274 | 0.98991 | 0.99739 |
| 0.1 | 0.84096 | 0.97644 | 0.55163 | 0.50990 | 0.36364 | 0.93414 | 0.91129 | 0.97558 |

Table 6 represents the targets of KAM when $\varepsilon = 0.0001$. According to the table, A should seriously improve the numbers of flights and passengers, even if the input factors are increased and the value of cargo is decreased. In other words, the inefficiency of A is due to the small numbers of flights and passengers, according to the data of other airports. These assessments are approximations and do not mean that it is impossible to increase the number of flights of airport A without increasing the values of its input factors, but at the same time, express the way of increasing the efficiency of A.

Table 6: The $10^{-4}$-KAM targets where $s_4^- \leq t\varepsilon x_{l4}$.

| Airport | Area | Apron | Terminal | Runway | Flights | Passengers | Cargo |
|---|---|---|---|---|---|---|---|
| A | 1,200.00 | 304,211.34 | 45,604.56 | 353,615.42 | 30,723.69 | 4,032,130.88 | 74,176.58 |
| B | 503.05 | 213,745.69 | 38,780.01 | 348,120.00 | 46,879.94 | 4,783,523.97 | 19,051.13 |
| C | 799.62 | 41,007.10 | 11,801.18 | 269,955.00 | 15,616.31 | 1,040,185.66 | 1,589.09 |
| D | 1,031.88 | 112,475.25 | 21,052.11 | 395,769.57 | 39,867.01 | 1,769,654.15 | 5,046.57 |
| E | 997.61 | 30,003.00 | 8,000.80 | 192,330.00 | 4,949.97 | 430,744.19 | 1,573.84 |
| F | 478.02 | 63,006.30 | 23,000.63 | 389,115.00 | 41,087.83 | 2,165,624.63 | 5,415.77 |
| G | 481.01 | 47,214.72 | 9,300.93 | 268,995.00 | 19,010.48 | 971,389.00 | 3,827.39 |
| H | 1,346.13 | 503,280.12 | 76,370.76 | 421,340.40 | 129,140.08 | 11,708,951.44 | 39,571.86 |

The same illustration can be discussed for other airports as well. In ad-dition, instead of SBM approach, the average measurement approach can also be used, that is, $W_j^- = 1/\text{ave}\{x_{ij}: i = 1,2,\ldots,8\}$, for $j = 1,2,3,4$, $W_k^- = 1/\text{ave}\{y_{ik}: i = 1,2,\ldots,8\}$, for $k = 1,2,3$, and so on for any other interested approaches introduced by expert judgment.

## 5.0. Conclusion

The meaning of technical efficiency, efficiency, effectiveness and productivity are discussed. Productivity is a combination of effectiveness and efficiency. In the effectiveness measurement, the factors of each firm are compared with desired goals of the firms, and in efficiency measurement the firms' performances are compared to each other. The provided technical efficiency scores should not be used to rank firms; they have unfortunately been used in the literature of operations research for the last four decades. The provided scores by KAM can relatively be meaningful according to the epsilon parameter and can be used to measure the efficiency scores of firms which lets us rank and benchmark the firms logically.